\documentclass[a4paper]{article}

\usepackage[margin=1in]{geometry} 

\usepackage[T1]{fontenc}
\newcommand{\changefont}[3]{
\fontfamily{#1} \fontseries{#2} \fontshape{#3} \selectfont}

\changefont{ptm}{m}{n}

\usepackage{setspace} 
\usepackage{graphicx}

\usepackage{amsfonts}
\usepackage{amsmath}
\usepackage{amssymb}
\usepackage{graphics}
\usepackage{mathrsfs}
\usepackage{color}
\newtheorem{remark}{Remark}[section]

\newtheorem{theorem}{Theorem}[section]

\newtheorem{lemma}{Lemma}[section]

\newtheorem{definition}{Definition}[section]

\long\def\symbolfootnote[#1]#2{\begingroup%
\def\thefootnote{\fnsymbol{footnote}}\footnote[#1]{#2}\endgroup} 

\begin{document}

\begin{center}
\Large \textbf{Unpredictability in Perturbed Quasilinear Systems with Regular Moments of Impulses}
\end{center}

\begin{center}
\normalsize \textbf{Mehmet Onur Fen$^{1,}\symbolfootnote[1]{Corresponding Author. E-mail: monur.fen@gmail.com}$, Fatma Tokmak Fen$^2$} \\
\vspace{0.2cm}
\textit{\textbf{$^1$Department of Mathematics, TED University, 06420 Ankara, Turkey}} \\

\vspace{0.1cm}
\textit{\textbf{$^2$Department of Mathematics, Gazi University, 06560 Ankara, Turkey}} \\
\vspace{0.1cm}
\end{center}

\vspace{0.3cm}

\begin{center}
\textbf{Abstract} 
\end{center}

It is rigorously proved that quasilinear impulsive systems possess unpredictable solutions when a perturbation generated by an unpredictable sequence is applied. The existence, uniqueness, as well as asymptotic stability of such solutions are demonstrated. The system under consideration is with regular moments of impulses, and for that reason a novel definition for unpredictable functions with regular discontinuity moments is provided. To show the existence of an unpredictable solution a Gronwall type inequality for piecewise continuous functions is utilized. The theoretical results are supported with an illustrative example. 

\vspace{-0.2cm}

\noindent\ignorespaces

\vspace{0.3cm}
 
\noindent\ignorespaces \textbf{Keywords:} Discontinuous unpredictable solution; Unpredictable sequence; Impulsive system; Asymptotic stability

\vspace{0.2cm}

\noindent\ignorespaces \textbf{MSC Classification:} 34A37, 34C10, 34C60 

\vspace{0.6cm}


\section{Introduction} \label{intro}

Impulsive differential equations are useful in modeling real world phenomena in which sudden changes occur. This type of differential equations have wide range of applications, for instance, in the fields of thermodynamics, neural networks, secure communication, economics, population dynamics, and control theory \cite{Haddad20}-\cite{Stamova16}. One of the basic and important problems in impulsive differential equations is the investigation of oscillations. Bounded, periodic, quasi-periodic, almost periodic, and anti-periodic oscillations in systems with impulses as well as their stability were investigated in the studies \cite{Li15}-\cite{Akh1}. 

The foundation for unpredictable oscillations was laid by Akhmet and Fen in the papers \cite{Fen16}-\cite{Fen18}. An unpredictable trajectory is a special type of Poisson stable one and its existence implies Poincar\'e chaos in the corresponding quasi-minimal set \cite{Fen16}. Differently from chaos in the senses of both Li-Yorke \cite{Li75} and Devaney \cite{Devaney87}, the definition of Poincar\'e chaos is based on a single motion, which is unpredictable. The logistic and H\'enon maps are examples of discrete dynamics comprising unpredictable orbits. The presence of continuous unpredictable oscillations in various types of differential equations was demonstrated in the studies \cite{Fen17a}-\cite{Fen18}. Moreover, the reader is referred to the papers \cite{Miller19}-\cite{Thakur21} for generalizations of Poincar\'e chaos and unpredictable points to topological spaces.
 
In the present study we investigate unpredictable solutions of the impulsive system
\begin{eqnarray}
\begin{array}{l} \label{mainimpulsive}
\displaystyle \frac{dx}{dt} = A x + f(t,x) + g(t),  ~t \neq  \theta_k, \\
\Delta x |_{t = \theta_k} = B x(\theta_k) + h(x(\theta_k)),
\end{array}
\end{eqnarray}
where the functions $f:\mathbb R\times \mathbb R^m\to \mathbb R^m$ and $h:\mathbb R^m \to \mathbb R^m$ are continuous in all their arguments, the function $f(t,x)$ is $\omega-$periodic in $t$ for some positive number $\omega$, i.e. $f(t+\omega,x)=f(t,x)$ for all $t\in\mathbb R$ and $x \in \mathbb R^m$, $A \in \mathbb R^{m \times m}$ and $B\in \mathbb R^{m \times m}$ are constant matrices, there exists a natural number $p$ such that the strictly increasing sequence $\left\{\theta_k\right\}_{k\in \mathbb Z}$ of impulse moments satisfies  
\begin{eqnarray} \label{regularitycond}
\theta_{k+p}=\theta_k+\omega
\end{eqnarray}
for all $k\in\mathbb Z$, $\Delta x |_{t = \theta_k} = x(\theta_k+) - x(\theta_k)$, and $x(\theta_k + )=\displaystyle \lim_{t\to \theta_k^+} x(t)$. The piecewise constant function $g : \mathbb R \to \mathbb R^m$ is defined by the equation
\begin{eqnarray} \label{discontineqn}
g(t) = \sigma_k
\end{eqnarray}
for $t \in (\theta_{kp}, \theta_{(k+1)p}]$, $k \in \mathbb Z$, in which $\left\{\sigma_k\right\}_{k \in \mathbb Z}$ is a bounded sequence in $\mathbb R^m$ with $\displaystyle\sup_{k \in \mathbb Z}  \left\|\sigma_k \right\| \leq M_{\sigma}$ for some positive number $M_{\sigma}$. It is worth noting that (\ref{mainimpulsive}) is obtained by applying the perturbation function $g(t)$ to the impulsive system 
\begin{eqnarray}
\begin{array}{l} \label{mainimpulsive2}
\displaystyle \frac{dx}{dt} = A x + f(t,x), ~t \neq  \theta_k, \\
\Delta x |_{t = \theta_k} = B x(\theta_k) + h(x(\theta_k)).
\end{array}
\end{eqnarray}
Our purpose is to prove that if the sequence $\left\{\sigma_k\right\}_{k \in \mathbb Z}$ is unpredictable, then the impulsive system (\ref{mainimpulsive}) possesses a unique asymptotically stable unpredictable solution under the sufficient conditions given in Section \ref{sec_prelim}. In other words, by considering the function $g(t)$ as an input applied to (\ref{mainimpulsive2}), we will show that an unpredictable input generates an unpredictable output. The definitions of an unpredictable sequence and piecewise continuous unpredictable function are given in Section \ref{sec3}.

Unpredictable solutions of linear impulsive systems were investigated in paper \cite{Fen21}, where the model admits unpredictable impulse moments. In contrast to paper \cite{Fen21}, in this study we consider an impulsive system which comprises nonlinear terms in both the differential and impulse equations. Another difference of our study is the usage of regular impulse moments in the model, i.e., the equation (\ref{regularitycond}) is satisfied. The concept of B-topology was utilized in \cite{Fen21} to define piecewise continuous unpredictable solutions since the unpredictability of the impulse moments causes the discontinuities of the solutions not to coincide under shifts of the time argument. However, in our case, the discontinuities of the solutions coincide under shifts of the time argument by integer multiples of $\omega$ because of the regularity of impulse moments. Accordingly, we make use of a novel definition for piecewise continuous unpredictable solutions with regular discontinuity moments regardless of the B-topology. Moreover, owing to the presence of nonlinear terms in the differential as well as impulse equations in (\ref{mainimpulsive}), we take advantage of a Gronwall type inequality for piecewise continuous functions to obtain the main result. Additionally, for the same reason, a more complicated proof technique is required compared to \cite{Fen21}. The results of the present study is also different compared to \cite{new} since impulsive systems with a discrete unpredictable set of discontinuity moments is taken into account in \cite{new}, whereas in our paper the system is with regular moments of impulses and the perturbation is obtained via an unpredictable sequence.

The remaining parts of the paper are organized as follows. In Section \ref{sec_prelim}, some preliminary results concerning the bounded solution of system (\ref{mainimpulsive}) are mentioned, and the assumptions which guarantee the existence, uniqueness and asymptotic stability of an unpredictable solution are given. Section \ref{sec3}, on the other hand, is concerned with the main result of the present study. In that section, the novel definition of an unpredictable function with regular discontinuity moments is provided, and the existence of a unique unpredictable solution of (\ref{mainimpulsive}), which is asymptotically stable, is proved. The theoretical results are supported with an illustrative example in Section \ref{examplessec}. In the example, the logistic map is used as the source of the unpredictability. Finally, Section \ref{secconc} is devoted to some concluding remarks.

\section{Preliminaries} \label{sec_prelim}
 
Throughout the paper we will make use of the Euclidean norm for vectors and the spectral norm for square matrices.

The following assumptions are required.
\begin{enumerate}
\item[\bf (A1)] The matrices $A$ and $B$ commute, and $\det(I+B)\neq 0,$ where $I$ is the $m \times m$ identity matrix;

\item[\bf (A2)] The eigenvalues of the matrix $A+\displaystyle \frac{p}{\omega} \log (I+B)$ have negative real parts;

\item[\bf (A3)] There exist positive numbers $M_f$ and $M_h$ such that $$\displaystyle \sup_{t\in \mathbb R, x\in \mathbb R^m} \left\|f(t,x)\right\| \leq M_f$$ and $$\displaystyle \sup_{x\in \mathbb R^m} \left\|h(x)\right\|\le M_h;$$

\item[\bf (A4)] There exist positive numbers $L_f$ and $L_h$ such that $$\left\|f(t,x_1)-f(t,x_2)\right\| \leq L_f\left\|x_1-x_2\right\|$$ for each $t\in\mathbb R,$ $x_1, x_2 \in \mathbb R^{m}$ and $$\left\|h(x_1)-h(x_2)\right\| \leq L_h \left\|x_1-x_2\right\|$$ for each $x_1, x_2 \in \mathbb R^{m}$. 
\end{enumerate} 

In the remaining parts of the paper we denote by $i(J)$ the number of the terms of the sequence $\left\{\theta_k\right\}_{k\in \mathbb Z}$ of impulse moments which belong to an interval $J$. It can be confirmed that 
\begin{eqnarray} \label{numberterms}
\displaystyle i((a,b))\leq p+\frac{p}{\omega}(b-a),
\end{eqnarray}
for each real numbers $a$ and $b$ with $b>a$.
Let $U(t,s)$ be the matriciant of the linear homogeneous impulsive system
\begin{eqnarray} \label{linearimpulsive}
\begin{array}{l}
\displaystyle \frac{dx}{dt} = A x,   ~t \neq  \theta_k, \\
\Delta x |_{t = \theta_k} = B x (\theta_k).
\end{array}
\end{eqnarray}
The equation $U(t,s)=e^{A(t-s)}(I+B)^{i([s,t))}$ holds for $t>s$ and $U(s,s)=I$ provided that the assumption $(A1)$ is valid. Additionally,  there exist positive numbers $N$ and $\lambda$ such that 
\begin{eqnarray*} \label{transitionmatrixineq}
\left\|U(t,s)\right\| \leq Ne^{-\lambda (t-s)}
\end{eqnarray*} 
for $t \geq s$ under the assumptions $(A1)$ and $(A2)$ \cite{Samolienko95,Akh1}.

The following assumptions are also needed. 
\begin{enumerate}
\item[\bf (A5)] $\displaystyle N \left( \frac{L_f}{\lambda} + \frac{p L_h}{1-e^{- \lambda \omega}}  \right) <1;$

\item[\bf (A6)] $\displaystyle  N L_f + \frac{p}{\omega} \ln (1+N L_h) < \lambda ;$

\item[\bf (A7)] $L_h \big\| \left(I+B\right)^{-1} \big\|<1.$
\end{enumerate}

According to the results of the book \cite{Samolienko95}, if the assumptions $(A1)-(A5)$ are valid, then the impulsive system (\ref{mainimpulsive}) admits a unique solution $\phi(t)$ bounded on the whole real axis which satisfies the relation
\begin{eqnarray} \label{integrosumeqn} 
\phi(t)=\displaystyle \int_{-\infty}^t U(t,s) \left[ f(s,\phi(s)) + g(s)\right] ds + \sum_{-\infty< \theta_k<t} U(t,\theta_k+) h (\phi(\theta_k)).
\end{eqnarray} 
By means of the estimates
$$ \Big\|\displaystyle \int_{-\infty}^t U(t,s) \left[ f(s,\phi(s)) + g(s)\right] ds \Big\| \leq \displaystyle \frac{N (M_f+M_{\sigma})}{\lambda}$$
and
$$ \Big\| \sum_{-\infty< \theta_k<t} U(t,\theta_k+) h (\phi(\theta_k)) \Big\| \leq \frac{p N M_h}{1-e^{-\lambda \omega}},$$
one can confirm that
$$\displaystyle \sup_{t \in \mathbb R} \left\| \phi (t)\right\| \leq M_{\phi},$$ where 
\begin{eqnarray} \label{Mphi} 
M_{\phi}= \displaystyle \frac{N (M_f+M_{\sigma})}{\lambda} + \frac{p N M_h}{1-e^{-\lambda \omega}}.
\end{eqnarray}

\section{Existence of an Unpredictable Solution} \label{sec3}

The definition of an unpredictable sequence is as follows.

\begin{definition} (\cite{Fen18}) \label{unpseqquasi}
A bounded sequence $\left\{\sigma_k\right\}_{k \in \mathbb Z}$ is called unpredictable if there exist a positive number $\delta_{0}$ (the unpredictability constant) and sequences $\left\{\zeta_n\right\}_{n\in\mathbb N},$ $\left\{\eta_n\right\}_{n\in\mathbb N}$ of positive integers both of which diverge to infinity such that $\left\|\sigma_{k+\zeta_n} - \sigma_k \right\|\to 0$ as $n \to \infty$ for each $k$ in bounded intervals of integers and $ \left\|\sigma_{\zeta_n + \eta_n} - \sigma_{\eta_n}\right\| \geq \delta_{0}$ for each $n\in\mathbb N.$
\end{definition}

The description of a discontinuous unpredictable function with regularly changing discontinuity moments is provided in the next definition.

\begin{definition}\label{unpfuncimp}
Suppose that there exist a natural number $p$ and a positive number $\omega$ such that the strictly increasing sequence $\left\{\theta_k\right\}_{k \in \mathbb Z}$ of real numbers satisfies the equation $\theta_{k+p}=\theta_k+\omega$ for each $k \in \mathbb Z$. A piecewise continuous and bounded function $\varphi : \mathbb  R  \to \mathbb R^m$ with the set of possible discontinuity points $\left\{\theta_k\right\}_{k \in \mathbb Z}$ satisfying $\displaystyle \lim_{t \to \theta_{k}^-} \varphi(t)=\varphi(\theta_{k})$, $k \in \mathbb{Z}$, is called unpredictable if there exist positive numbers $\epsilon_{0}$ (the unpredictability constant), $r$ and sequences $\left\{\mu_{n}\right\}_{n\in\mathbb N}$, $\left\{\tau_{n}\right\}_{n\in\mathbb N}$ of real numbers both of which diverge to infinity such that:
\begin{itemize}	
\item[(i)] for every positive number $\epsilon$ there exists a positive number $\delta$ such that $\left\|\varphi(s_1)-\varphi(s_2) \right\|<\epsilon $ whenever the points $s_1$ and $s_2$ belong to the same interval of continuity and $\left|s_1-s_2\right| <\delta$;
\item[(ii)] $\left\|\varphi(t+\mu_n) - \varphi(t) \right\| \to 0 $ as $n \to \infty$ uniformly on compact subsets of $\mathbb R$;
\item[(iii)] $\|\varphi(t+\mu_{n})-\varphi(t)\| \geq \epsilon_{0}$ for each $t\in [\tau_{n}-r, \tau_{n}+r]$ and $n \in \mathbb N$.
	\end{itemize}	
\end{definition}



The following analogue of Gronwall's inequality for piecewise continuous functions is required.
\begin{lemma} (\cite{Bainov92}) \label{gronwalllemma}
Suppose that for $t_1 \leq t \leq t_2$ the following inequality holds:
$$u(t) \leq a(t) + \displaystyle \int_{t_1}^t b(s) u(s) ds + \sum_{t_1 < \theta_k < t} \beta_k u(\theta_k),$$
where $\beta_k \geq 0$, $k\in \mathbb N$, are constants and $u(t):\mathbb R \to \mathbb R$, $a(t):\mathbb R \to \mathbb R$, $b(t):\mathbb R \to [0,\infty)$ are piecewise continuous functions that have first kind discontinuities at the points $\theta_k$, $k\in\mathbb Z$, only and are left continuous at each $\theta_k$. Then, for $t_1 \leq t \leq t_2$,
\begin{eqnarray*}
u(t) \leq a(t) + \displaystyle \int_{t_1}^t a(s) b(s) \prod_{s < \theta_k < t} (1+ \beta_k) e^{\int_s^t b(\tau) d\tau} ds + \sum_{t_1 < \theta_k < t} a(\theta_k) \beta_k \prod_{\theta_k < \theta_j < t} (1+\beta_j) e^{\int_{\theta_k}^t b(\tau) d\tau}.
\end{eqnarray*}
\end{lemma}
 
We make use of Lemma \ref{gronwalllemma} in the proof of the next assertion.
\begin{lemma} \label{impulselemma}
Suppose that there is a sequence $\left\{\zeta_n\right\}_{n\in\mathbb N}$ of positive integers which diverges to infinity such that $\left\|\sigma_{k+\zeta_n} - \sigma_{k} \right\|\to 0$ as $n \to \infty$ for each $k$ in bounded intervals of integers. Then, under the assumptions $(A1)-(A6)$, the bounded solution $\phi(t)$ of impulsive system (\ref{mainimpulsive}) has the property that $\left\|\phi(t+\omega \zeta_n)-\phi(t)\right\| \to 0$ as $n \to \infty$ uniformly on compact subsets of $\mathbb R$.
\end{lemma} 
 
\noindent {\bf Proof.}  Let us fix a positive number $\epsilon$ and a compact subset $\mathcal{C}$ of the real axis. There exist integers $\alpha$ and $\beta$ with $\beta > \alpha$ such that $\mathcal{C} \subseteq \left[ \theta_{\alpha p}, \theta_{\beta p} \right]$.
We define the numbers 
\begin{eqnarray} \label{proofnum1}
K_1 = \displaystyle \left[ \frac{2N (M_f + M_{\sigma})}{\lambda} + \frac{2 p N M_h}{1- e^{-\lambda \omega}}  \right] (1+ N L_{h})^p 
\end{eqnarray}
and 
\begin{eqnarray} \label{proofnum2}
K_2 = \frac{N}{\lambda} + \frac{N^2 L_{f} (1+ N L_{h})^p}{\lambda \gamma} + \frac{N^2 p L_{h} (1+N L_{h})^p e^{\gamma \omega}}{\lambda \left(1-e^{-\gamma \omega}\right)},
\end{eqnarray}
where
\begin{eqnarray} \label{proofnum3} 
\gamma = \lambda - N L_{f} - \frac{p}{\omega} \ln (1+N L_{h}).  
\end{eqnarray}  
It is worth noting that the number $\gamma$ is positive according to assumption $(A6)$.  
  
Next, we take a positive number $\rho_0$ such that 
\begin{eqnarray} \label{proofineq1} 
\rho_0 \leq \frac{1}{K_1 + K_2}
\end{eqnarray}  
and a positive integer $j$ satisfying 
\begin{eqnarray} \label{proofineq2} 
j \geq \displaystyle \frac{1}{\gamma \omega} \ln \left( \frac{1}{\rho_0 \epsilon} \right).
\end{eqnarray}  
Since $\left\|\sigma_{k+\zeta_n} - \sigma_{k} \right\|\to 0$ as $n \to \infty$ for each $k$ in bounded intervals of integers, there exists a natural number $n_0$ such that for $n \geq n_0$ we have
\begin{eqnarray} \label{proofineq3} 
\left\|\sigma_{k+ \zeta_n} - \sigma_k\right\| < \rho_0 \epsilon 
\end{eqnarray} 
for each $k=\alpha-j, \ \alpha-j+1, \ \cdots, \ \beta-1$. 

Let us fix a natural number $n \geq n_0$. The inequality
\begin{eqnarray} \label{proofineq4} 
\left\|g(t+ \omega \zeta_n) - g(t)\right\| < \rho_0 \epsilon, \ \ \theta_{(\alpha-j)p} < t \leq \theta_{\beta p}. 
\end{eqnarray}   
holds in compliance with (\ref{proofineq3}).

It can be verified with the aid of (\ref{integrosumeqn}) that 
\begin{eqnarray*} 
\phi(t + \omega \zeta_n) - \phi(t) & = & \displaystyle \int_{- \infty}^{t} U(t,s) \big( f(s, \phi(s+\omega \zeta_n)) - f(s,\phi(s)) + g(s+\omega \zeta_n) - g(s) \big) \ ds \\
& +& \sum_{-\infty < \theta_k < t} U(t, \theta_k+) \big( h(\phi(\theta_k+\omega \zeta_n)) - h(\phi(\theta_k)) \big).
\end{eqnarray*}

Now, suppose that $t \in \left[\theta_{(\alpha-j)p}, \ \theta_{\beta p}\right]$.
The inequality (\ref{proofineq4}) yields
\begin{eqnarray*} 
\Big\| \displaystyle  \int^{t}_{\theta_{(\alpha-j)p}} U(t,s)  \big( g(s+\omega \zeta_n) -g(s) \big) \ ds \Big\|
\leq \frac{N \rho_0 \epsilon}{\lambda} \left(1-e^{-\lambda (t- \theta_{(\alpha-j)p} )}\right).
\end{eqnarray*}    
In addition to this, the estimates 
\begin{eqnarray*}  
&& \Big\| \displaystyle  \int_{- \infty}^{\theta_{(\alpha-j)p}} U(t,s) \big( f(s, \phi(s+\omega \zeta_n)) - f(s,\phi(s)) + g(s+\omega \zeta_n) - g(s) \big) \ ds \Big\| \\
&& \leq \frac{2 N (M_f +M_{\sigma})}{\lambda} \ e^{-\lambda (t-\theta_{(\alpha-j)p})}
\end{eqnarray*}  
and
\begin{eqnarray*}
\Big\| \displaystyle  \sum_{-\infty < \theta_k \leq \theta_{(\alpha-j)p}} U(t, \theta_k+) \big( h(\phi(\theta_k+\omega \zeta_n)) - h(\phi(\theta_k)) \big) \Big\| 
\leq \frac{2 p N M_{h}}{1-e^{-\lambda \omega}} e^{-\lambda (t-\theta_{(\alpha-j)p})}
\end{eqnarray*}   
are fulfilled. Thus, 
\begin{eqnarray*}
\left\|\phi(t + \omega \zeta_n) - \phi(t)\right\| & \leq & \left[ \frac{2N (M_f + M_{\sigma})}{\lambda} + \frac{2 p N M_h}{1- e^{-\lambda \omega}}  \right] e^{-\lambda (t-\theta_{(\alpha-j)p})}  + \frac{N \rho_0 \epsilon}{\lambda} \left(1-e^{-\lambda (t-\theta_{(\alpha-j)p})}\right) \\
& + & \int_{\theta_{(\alpha-j)p}}^{t} N L_{f} e^{- \lambda (t-s)} \left\|\phi(s + \omega \zeta_n) - \phi(s)\right\| ds \\
& + & \sum_{\theta_{(\alpha-j)p} < \theta_k < t} N L_{h} e^{-\lambda (t-\theta_k)} \left\|\phi(\theta_k + \omega \zeta_n) - \phi(\theta_k)\right\|.
\end{eqnarray*}  
Accordingly, we have 
\begin{eqnarray*} 
u(t) \leq c +  \frac{N \rho_0 \epsilon}{\lambda} e^{\lambda t} 
 +  \int_{\theta_{(\alpha-j)p}}^{t} N L_{f} u(s) ds + \sum_{\theta_{(\alpha-j)p} < \theta_k < t} N L_{h} u(\theta_k),
\end{eqnarray*}
where $$u(t) = e^{\lambda t} \left\|\phi(t + \omega \zeta_n) - \phi(t)\right\| $$ 
and
$$ 
c=\left[ \frac{2N (M_f + M_{\sigma})}{\lambda} + \frac{2 p N M_h}{1- e^{-\lambda \omega}} - \frac{N \rho_0 \epsilon}{\lambda} \right] e^{\lambda  \theta_{(\alpha-j)p}}.
$$
By means of Lemma \ref{gronwalllemma} we obtain
\begin{eqnarray*}  
u(t) & \leq & c (1+ N L_{h})^{i((\theta_{(\alpha-j)p}, t))} e^{N L_{f} (t- \theta_{(\alpha-j)p})} + \frac{N \rho_0 \epsilon}{\lambda} e^{\lambda t} \\
& + & \frac{N^2 L_{f} \rho_0 \epsilon}{\lambda} \int_{\theta_{(\alpha-j)p}}^t (1+N L_{h})^{i((s,t))} e^{\lambda s} e^{N L_{f} (t-s)} ds \\
& + & \frac{N^2 L_{h} \rho_0 \epsilon}{\lambda} \sum_{\theta_{(\alpha-j)p} < \theta_k < t} (1+N L_{h})^{i((\theta_k, t))} e^{\lambda \theta_k} e^{N L_{f} (t-\theta_k)}.
\end{eqnarray*} 
One can attain using (\ref{numberterms}) that
\begin{eqnarray*} 
\int_{\theta_{(\alpha-j)p}}^t (1+N L_{h})^{i((s,t))} e^{\lambda s} e^{N L_{f} (t-s)} ds 
\leq \frac{(1+N L_{h})^p}{\gamma} \ e^{\lambda t}  \left(1-e^{-\gamma (t-\theta_{(\alpha-j)p}) }\right)
\end{eqnarray*}
and
\begin{eqnarray*} 
\sum_{\theta_{(\alpha-j)p} < \theta_k < t} (1+N L_{h})^{i((\theta_k, t))} e^{\lambda \theta_k} e^{N L_{f} (t-\theta_k)}
\leq \frac{p (1+N L_{h})^p e^{\gamma \omega}}{1-e^{-\gamma \omega}} \ e^{\lambda t} \left(1-e^{-\gamma (t-\theta_{(\alpha-j-1)p}) }\right),
\end{eqnarray*} 
 where $\gamma$ is the number defined by (\ref{proofnum3}). Hence, we have
 \begin{eqnarray*} 
 u(t) & \leq & c (1+ N L_{h})^p e^{(\lambda - \gamma)(t-\theta_{(\alpha-j)p})} + \frac{N \rho_0 \epsilon}{\lambda} e^{\lambda t} \\
& + & \frac{N^2 L_{f} (1+N L_{h})^p \rho_0 \epsilon}{\lambda \gamma} \ e^{\lambda t}  \left(1-e^{-\gamma (t-\theta_{(\alpha-j)p}) }\right) \\
& + & \frac{N^2 p L_{h} (1+N L_{h})^p e^{\gamma \omega} \rho_0 \epsilon}{\lambda (1-e^{-\gamma \omega})} \ e^{\lambda t} \left(1-e^{-\gamma (t-\theta_{(\alpha-j-1)p}) }\right).
 \end{eqnarray*} 
The last inequality implies that
\begin{eqnarray*} 
\left\|\phi(t + \omega \zeta_n) - \phi(t)\right\| < K_1 e^{-\gamma (t-\theta_{(\alpha-j)p})} + K_2 \rho_0 \epsilon, \ \    \theta_{(\alpha-j)p} \leq t \leq \theta_{\beta p}, 
\end{eqnarray*} 
in which $K_1$ and $K_2$ are respectively defined by (\ref{proofnum1}) and (\ref{proofnum2}).
 
If $\theta_{\alpha p} \leq t \leq \theta_{\beta p}$, then it can be verified using (\ref{proofineq2}) that $e^{-\gamma (t-\theta_{(\alpha-j)p})} \leq \rho_0 \epsilon$. Since $\mathcal{C} \subseteq \left[\theta_{\alpha p}, \theta_{\beta p}\right]$, in accordance with (\ref{proofineq1}), we have
$$\left\|\phi(t + \omega \zeta_n) - \phi(t)\right\| < (K_1 + K_2) \rho_0 \epsilon \leq \epsilon$$ for $t \in \mathcal{C}$. Consequently, $\left\|\phi(t + \omega \zeta_n) - \phi(t)\right\| \to 0$ as $n \to \infty$ uniformly on compact subsets of $\mathbb R$. $\square$

The main result of the present paper is given in the next theorem.

\begin{theorem} \label{impulsethm}
Suppose that the assumptions $(A1)-(A7)$ hold. If the sequence $\left\{\sigma_k\right\}_{k \in \mathbb Z}$ is unpredictable, then the impulsive system (\ref{mainimpulsive}) possesses a unique asymptotically stable unpredictable solution. 
\end{theorem}

\noindent {\bf Proof.} In the proof we will show that the bounded solution $\phi(t)$ of (\ref{mainimpulsive}) is unpredictable in the sense of Definition \ref{unpfuncimp}. According to the results of the book \cite{Samolienko95}, the bounded solution $\phi(t)$ is asymptotically stable provided that the assumptions $(A1)-(A6)$ are valid. If $t \in (\theta_k, \theta_{k+1})$ for some fixed integer $k$, then $\left\|\displaystyle \frac{d\phi}{dt}\right\| \leq \left\|A\right\| M_{\phi} + M_f + M_{\sigma}$, where $M_{\phi}$ is defined by (\ref{Mphi}). Therefore, Definition \ref{unpfuncimp}, (i) is satisfied for $\phi(t)$.

Since $\left\{\sigma_k\right\}_{k \in \mathbb Z}$ is an unpredictable sequence, there exist a positive number $\delta_0$ and sequences $\left\{\zeta_n\right\}_{n \in \mathbb N}$, $\left\{\eta_n\right\}_{n \in \mathbb N}$ of positive integers both of which diverge to infinity such that $\left\|\sigma_{k + \zeta_n} - \sigma_k \right\| \to 0$ as $n \to \infty$ for each $k$ in bounded intervals of integers and 
\begin{eqnarray} \label{thmproof1}
\left\|\sigma_{\zeta_n + \eta_n} - \sigma_{\eta_n}\right\| \geq \delta_0
\end{eqnarray} 
for each $n \in \mathbb N$.

Let us denote $$\mu_n= \omega \zeta_n$$ for each $n \in \mathbb N$. According to Lemma \ref{impulselemma}, $\left\|\phi(t+\mu_n) - \phi(t) \right\| \to 0 $ as $n \to \infty$ uniformly on compact subsets of $\mathbb R$. Since $\left\{\zeta_n\right\}_{n \in \mathbb N}$ diverges to infinity, the same is true for the sequence $\left\{\mu_n\right\}_{n \in \mathbb N}$.

In the remaining parts of the proof, we will show the existence of positive numbers $\epsilon_0$, $r$ and a sequence $\left\{\tau_n\right\}_{n \in \mathbb N}$, which diverges to infinity, such that
\begin{eqnarray} \label{thmproof2}
\left\| \phi(t + \mu_n) - \phi(t) \right\| \geq \epsilon_0
\end{eqnarray} 
for each $t \in [\tau_n-r, \tau_n+r]$ and $n \in \mathbb N$.

Let $\sigma_k=\left(\sigma_k^1,\sigma_k^2, \ldots, \sigma_k^m\right) \in \mathbb R^m$ for each $k \in \mathbb N$, and fix a number $n \in \mathbb N$.
The inequality (\ref{thmproof1}) yields 
$$
\big| \sigma_{\zeta_n + \eta_n}^{j_0} - \sigma_{\eta_n}^{j_0} \big| \geq \displaystyle \frac{\delta_0}{\sqrt{m}}
$$
for some integer $j_0$ with $1 \leq j_0 \leq m$.
Because the equation
$
g(t + \mu_n) - g(t) = \sigma_{\zeta_n + \eta_n} - \sigma_{\eta_n}
$
holds for $\theta_{\eta_n p}<t\leq\theta_{(\eta_n + 1)p}$, we have
\begin{eqnarray} \label{thmproofineq1}
\Big\| \displaystyle \int_{\theta_{\eta_n p}}^{\theta_{(\eta_n + 1)p}}  \left(g(s + \mu_n) -g(s) \right) ds \Big\|
=  \omega \bigg( \displaystyle \sum_{j=1}^{m} \left( \sigma_{\zeta_n + \eta_n}^j - \sigma^j_{\eta_n} \right)^2 \bigg)^{1/2}
\geq \omega \big| \sigma_{\zeta_n + \eta_n}^{j_0} - \sigma^{j_0}_{\eta_n} \big|
\geq \displaystyle \frac{\omega \delta_0}{\sqrt{m}}.
\end{eqnarray}
Making use of the equations
\begin{eqnarray*}
\phi(\theta_{(\eta_n+1)p}) & = & \phi(\theta_{\eta_n p}) + \displaystyle \int_{\theta_{\eta_n p}}^{\theta_{(\eta_n+1) p}} \left( A \phi(s) +f(s, \phi(s)) + g(s) \right) ds \\
& + &  \displaystyle \sum_{\theta_{\eta_n p} \leq \theta_k < \theta_{(\eta_n+1)p}} \left(B \phi(\theta_k) + h(\phi(\theta_k))\right)
\end{eqnarray*}
and
\begin{eqnarray*}
\phi(\theta_{(\eta_n+1)p} + \mu_n) & = & \phi(\theta_{\eta_n p} + \mu_n) + \displaystyle \int_{\theta_{\eta_n p}}^{\theta_{(\eta_n+1) p}} \left( A \phi(s+ \mu_n) +f(s, \phi(s+ \mu_n)) + g(s+ \mu_n) \right) ds \\
& + &  \displaystyle \sum_{\theta_{\eta_n p} \leq \theta_k < \theta_{(\eta_n+1)p}} \left(B \phi(\theta_k + \mu_n) + h(\phi(\theta_k+ \mu_n))\right)
\end{eqnarray*}
one can obtain that
\begin{eqnarray*}
 \left\| \phi(\theta_{(\eta_n+1)p} + \mu_n) - \phi(\theta_{(\eta_n+1)p}) \right\| 
& \geq & \Big\| \displaystyle \int_{\theta_{\eta_n p}}^{\theta_{(\eta_n + 1)p}}  \left(g(s + \mu_n) -g(s) \right) ds \Big\|
 - \left\| \phi(\theta_{\eta_n p} + \mu_n) - \phi(\theta_{\eta_n p}) \right\| \\
 & - & \displaystyle \int_{\theta_{\eta_n p}}^{\theta_{(\eta_n + 1)p}} \left(\left\|A\right\| + L_f \right) \left\|\phi(s+ \mu_n) - \phi(s) \right\| ds \\
 & - & \displaystyle \sum_{\theta_{\eta_n p} \leq \theta_k < \theta_{(\eta_n+1)p}} \left( \left\|B\right\| + L_h\right) \left\|\phi( \theta_k+ \mu_n) - \phi(\theta_k) \right\|. 
\end{eqnarray*}
According to (\ref{thmproofineq1}), we have that 
\begin{eqnarray*} 
\left\| \phi(\theta_{(\eta_n+1)p} + \mu_n) - \phi(\theta_{(\eta_n+1)p}) \right\|  
& \geq & \displaystyle \frac{\omega \delta_0}{\sqrt{m}}  - \left\| \phi(\theta_{\eta_n p} + \mu_n) - \phi(\theta_{\eta_n p}) \right\| \\
& - & \omega \left(\left\|A\right\| + L_f\right) \sup_{t \in [\theta_{\eta_n p}, \theta_{(\eta_n+1) p}]} \left\|\phi(t+\mu_n) - \phi(t)\right\|  \\
& - & p \left(\left\|B\right\| + L_h\right) \sup_{t \in [\theta_{\eta_n p}, \theta_{(\eta_n+1) p}]} \left\|\phi(t+\mu_n) - \phi(t)\right\|.
\end{eqnarray*} 
Thus,
\begin{eqnarray*}  
\sup_{t \in [\theta_{\eta_n p}, \theta_{(\eta_n+1) p}]} \left\|\phi(t+\mu_n) - \phi(t)\right\| \geq H_0,
\end{eqnarray*}  
where
$$
H_0= \displaystyle \frac{\omega \delta_0}{\left[2 + \omega \left(\left\|A\right\| + L_f\right) +  p \left(\left\|B\right\| + L_h\right) \right] \sqrt{m}}.
$$ 
Let us denote
\begin{eqnarray} \label{unpconst}
\epsilon_0= \frac{H_0}{2} \displaystyle \min\left\{ 1, \frac{\left( 1- L_h \left\|(I+B)^{-1}\right\|\right)}{2 \left\|(I+B)^{-1}\right\|}, \frac{1}{2\left(\left\|I+B\right\|+L_h\right)}\right\}
\end{eqnarray}
and
$$
r= \displaystyle \min\left\{ \frac{\overline{\theta}}{3},  \frac{H_0}{4R_0}, \frac{H_0 \left(1-L_h \left\|(I+B)^{-1}\right\|\right)}{8 R_0 \left\|(I+B)^{-1}\right\|}, \frac{H_0}{8R_0\left(\left\|I+B\right\|+L_h\right)} \right\}, 
$$ 
in which $$R_0=\left\|A\right\| M_{\phi} + M_f + M_{\sigma}$$ and $$\overline{\theta} = \displaystyle \min_{1 \leq k \leq p} \left(\theta_{k+1} - \theta_k\right).$$ It is worth noting that a closed interval of length $2r$ can contain at most one term of the sequence $\left\{\theta_k\right\}_{k \in \mathbb Z}$ of impulse moments since $\overline{\theta} \geq 3r$. 
 
In order to show the existence of a sequence $\left\{\tau_n\right\}_{n \in \mathbb N}$, which diverges to infinity, such that (\ref{thmproof2}) is valid for each $t \in [\tau_n-r, \tau_n+r]$ and $n \in \mathbb N$, we take into account the following three cases.

\noindent {\bf Case I.} First of all, we consider the case
\begin{eqnarray*}
\sup_{t \in [\theta_{\eta_n p}, \theta_{(\eta_n+1) p}]} \left\|\phi(t+\mu_n) - \phi(t)\right\| = \left\|\phi(\tau_n + \mu_n) - \phi(\tau_n)\right\|,
\end{eqnarray*} 
where $\tau_n \in (\theta_{k_0}, \theta_{k_0+1})$ for some integer $k_0$, which depends on $n$, with $\eta_n p \leq k_0 \leq (\eta_n + 1)p-1$.
For $t \in (\theta_{k_0}, \theta_{k_0+1}]$, the equation
\begin{eqnarray} \label{thmproofeqn3}
\phi(t + \mu_n) - \phi(t) & =& \phi(\tau_n + \mu_n) - \phi(\tau_n) \nonumber \\
 & + & \displaystyle \int_{\tau_n}^t  \big[ A \left(\phi(s+\mu_n) - \phi(s) \right)  + f(s,\phi(s+\mu_n)) \\ && -  f(s,\phi(s)) + g(s+\mu_n) - g(s) \big] ds \nonumber
\end{eqnarray} 
is satisfied.
 
If $[\tau_n-r, \tau_n+r] \subset (\theta_{k_0}, \theta_{k_0+1}]$, then for each $t \in [\tau_n-r, \tau_n+r]$ we have
\begin{eqnarray} \label{thmproofeqn4}
\left\|\phi(t + \mu_n) - \phi(t)\right\| & \geq & \left\|\phi(\tau_n + \mu_n) - \phi(\tau_n)\right\|  \nonumber  \\
 & - & \displaystyle \Big| \int_{\tau_n}^t  \big\| A \left(\phi(s+\mu_n) - \phi(s) \right)  + f(s,\phi(s+\mu_n))  \nonumber \\ &&  -  f(s,\phi(s)) + g(s+\mu_n) - g(s) \big\| ds \Big| \\
 & \geq & H_0 - 2r R_0 \nonumber \\ 
 & \geq & \displaystyle \frac{H_0}{2} \geq \epsilon_0. \nonumber 
\end{eqnarray}  

Next, suppose that $\tau_n+r>\theta_{k_0+1}$. Then, in a similar way to (\ref{thmproofeqn4}), one can obtain by means of equation (\ref{thmproofeqn3}) that 
$
\left\|\phi(t + \mu_n) - \phi(t)\right\| \geq H_0/2
$
for $t \in [\tau_n -r, \theta_{k_0+1}]$. In particular, the inequality
$
\left\|\phi(\theta_{k_0+1} + \mu_n) - \phi(\theta_{k_0+1})\right\| \geq \displaystyle H_0/2
$ holds.
Benefiting from the assumption $(A7)$ we attain that
\begin{eqnarray*}
\left\|\phi(\theta_{k_0+1} + \mu_n  + ) - \phi(\theta_{k_0+1} + )\right\| & \geq &  \big\| (I+B) (\phi(\theta_{k_0+1} + \mu_n) - \phi(\theta_{k_0+1})) \big\| \\ & - &   \big\| h(\phi(\theta_{k_0+1} + \mu_n)) - h(\phi(\theta_{k_0+1})) \big\| \\
 &\geq & \displaystyle \frac{\left( 1- L_h \left\|(I+B)^{-1}\right\|\right)}{\left\|(I+B)^{-1}\right\|}   \big\| (\phi(\theta_{k_0+1} + \mu_n) - \phi(\theta_{k_0+1})) \big\| \\
 &\geq & \displaystyle \frac{H_0 \left( 1- L_h \left\|(I+B)^{-1}\right\|\right)}{2 \left\|(I+B)^{-1}\right\|}.
\end{eqnarray*}
Therefore, 
\begin{eqnarray*}
\left\|\phi(t + \mu_n) - \phi(t)\right\| & \geq & \left\|\phi(\theta_{k_0+1} + \mu_n + ) - \phi(\theta_{k_0+1} +)\right\| \\
& - & \displaystyle \int_{\theta_{k_0+1}}^t  \big\| A \left(\phi(s+\mu_n) - \phi(s) \right)  + f(s,\phi(s+\mu_n)) \\ && -  f(s,\phi(s)) + g(s+\mu_n) - g(s) \big\| ds \\
&>& \displaystyle \frac{H_0 \left( 1- L_h \left\|(I+B)^{-1}\right\|\right)}{2 \left\|(I+B)^{-1}\right\|} -2rR_0 \\
&\geq & \displaystyle \frac{H_0 \left( 1- L_h \left\|(I+B)^{-1}\right\|\right)}{4 \left\|(I+B)^{-1}\right\|}
\end{eqnarray*}
for $t \in (\theta_{k_0+1}, \tau_n + r]$.
For that reason, if $t\in [\tau_n-r, \tau_n+r]$, then
$$\left\|\phi(t+ \mu_n) - \phi(t)\right\| \geq \displaystyle \frac{H_0}{2} \min\left\{1, \frac{\left( 1- L_h \left\|(I+B)^{-1}\right\|\right)}{2 \left\|(I+B)^{-1}\right\|}\right\} \geq \epsilon_0. $$

Now, assume that $\tau_n-r \leq \theta_{k_0}$. Utilizing equation (\ref{thmproofeqn3}) one more time, it can be deduced for $t \in (\theta_{k_0}, \tau_n+r]$ that
$
\left\|\phi(t+ \mu_n) - \phi(t) \right\| \geq \displaystyle H_0/2,
$
and in particular,
$\left\|\phi(\theta_{k_0}+ \mu_n +) - \phi(\theta_{k_0}+) \right\| \geq \displaystyle H_0/2$. Accordingly, the inequality
$$\left\|\phi(\theta_{k_0}+ \mu_n) - \phi(\theta_{k_0}) \right\| \geq \displaystyle \frac{H_0}{2\left(\left\|I+B\right\|+L_h\right)}$$
is valid. Thus, if $\tau_n-r=\theta_{k_0}$, then for $t \in [\tau_n-r, \tau_n+r]$ we have
$$\left\|\phi(t+ \mu_n) - \phi(t) \right\| \geq \displaystyle \frac{H_0}{2}\min\left\{1,\frac{1}{\left\|I+B\right\|+L_h} \right\}\geq \epsilon_0.$$ On the other hand, if $\tau_n-r < \theta_{k_0}$, then for $t \in [\tau_n-r, \theta_{k_0}]$ one can obtain
\begin{eqnarray*}
\left\|\phi(t + \mu_n) - \phi(t)\right\| & \geq & \left\|\phi(\theta_{k_0} + \mu_n) - \phi(\theta_{k_0})\right\| \\
& - & \displaystyle \Big| \int_{\theta_{k_0}}^t  \big\| A \left(\phi(s+\mu_n) - \phi(s) \right) + f(s,\phi(s+\mu_n)) \\ && -  f(s,\phi(s)) + g(s+\mu_n) - g(s) \big\| ds \Big|\\
&>& \displaystyle \frac{H_0}{2\left(\left\|I+B\right\|+L_h\right)} - 2rR_0 \\
&\geq & \displaystyle \frac{H_0}{4\left(\left\|I+B\right\|+L_h\right)}.
\end{eqnarray*}
Hence, the inequality 
$$\left\|\phi(t+ \mu_n) - \phi(t) \right\| \geq \displaystyle \frac{H_0}{2}\min\left\{1,\frac{1}{2\left(\left\|I+B\right\|+L_h\right)} \right\}\geq \epsilon_0$$
holds for $t \in [\tau_n-r, \tau_n+r]$.
 
\noindent {\bf Case II.} Secondly, we suppose that
\begin{eqnarray*}  
\sup_{t \in [\theta_{\eta_n p}, \theta_{(\eta_n+1) p}]} \left\|\phi(t+\mu_n) - \phi(t)\right\| = \left\|\phi(\theta_{k_1} + \mu_n) - \phi(\theta_{k_1})\right\|
\end{eqnarray*}
for some integer $k_1$, which depends on $n$, satisfying $\eta_n p \leq k_1 \leq (\eta_n+1)p$. 
In this case,
\begin{eqnarray*}
\left\|\phi(\theta_{k_1} + \mu_n + ) - \phi(\theta_{k_1}+)\right\| & \geq & \left\|(I+B) \left(\phi(\theta_{k_1} + \mu_n) - \phi(\theta_{k_1})\right)\right\| \\
& - & \big\| h(\phi(\theta_{k_1} + \mu_n)) - h(\phi(\theta_{k_1})) \big\| \\
& \geq & \displaystyle \frac{\left( 1- L_h \left\|(I+B)^{-1}\right\|\right)}{\left\|(I+B)^{-1}\right\|}   \big\| (\phi(\theta_{k_1} + \mu_n) - \phi(\theta_{k_1})) \big\| \\
& \geq & \displaystyle \frac{H_0 \left( 1- L_h \left\|(I+B)^{-1}\right\|\right)}{\left\|(I+B)^{-1}\right\|}. 
\end{eqnarray*} 
 
Let us set $\tau_n=\theta_{k_1}$. For $t \in (\tau_n, \tau_n+r]$, we have
\begin{eqnarray*}
\left\|\phi(t+\mu_n) - \phi(t) \right\| & \geq &  \left\|\phi(\tau_n+\mu_n +) - \phi(\tau_n+) \right\| \\
& - & \displaystyle \int_{\tau_{n}}^t  \big\| A \left(\phi(s+\mu_n) - \phi(s) \right)  + f(s,\phi(s+\mu_n)) \\ && -  f(s,\phi(s)) + g(s+\mu_n) - g(s) \big\| ds  \\
& \geq & \displaystyle \frac{H_0 \left( 1- L_h \left\|(I+B)^{-1}\right\|\right)}{\left\|(I+B)^{-1}\right\|} -2rR_0 \\
  & \geq & \displaystyle \frac{3H_0 \left( 1- L_h \left\|(I+B)^{-1}\right\|\right)}{4\left\|(I+B)^{-1}\right\|}.
\end{eqnarray*}
For $t \in [\tau_n-r,\tau_n]$, on the other hand, one can attain that
\begin{eqnarray*}
\left\|\phi(t+\mu_n) - \phi(t) \right\| & \geq &  \left\|\phi(\tau_n+\mu_n) - \phi(\tau_n) \right\| \\
& - & \displaystyle \Big|\int_{\tau_{n}}^t  \big\| A \left(\phi(s+\mu_n) - \phi(s) \right)  + f(s,\phi(s+\mu_n)) \\ && -  f(s,\phi(s)) + g(s+\mu_n) - g(s) \big\| ds  \Big| \\
& \geq & H_0 -2rR_0 \\
 & \geq & \frac{H_0}{2}.
\end{eqnarray*} 
Thus, we have 
$$
\left\|\phi(t+\mu_n) - \phi(t) \right\| \geq \displaystyle \frac{H_0}{2} \min\left\{1,\displaystyle \frac{3 \left( 1- L_h \left\|(I+B)^{-1}\right\|\right)}{2\left\|(I+B)^{-1}\right\|}\right\} \geq \epsilon_0  
$$
for $t \in [\tau_n-r, \tau_n+r]$.

\noindent {\bf Case III.} Finally, we take into account the case 
\begin{eqnarray*}  
\sup_{t \in [\theta_{\eta_n p}, \theta_{(\eta_n+1) p}]} \left\|\phi(t+\mu_n) - \phi(t)\right\| = \left\|\phi(\theta_{k_2} + \mu_n+) - \phi(\theta_{k_2}+)\right\|
\end{eqnarray*}
for some integer $k_2$ depending on $n$ such that $\eta_n p \leq k_2 < (\eta_n+1)p$.
  
We set $\tau_n=\theta_{k_2}$. If $t \in (\tau_n, \tau_n+r]$, then
\begin{eqnarray*}
\left\|\phi(t+\mu_n) - \phi(t) \right\| & \geq &  \left\|\phi(\tau_n+\mu_n +) - \phi(\tau_n+) \right\| \\
& - & \displaystyle \int_{\tau_{n}}^t  \big\| A \left(\phi(s+\mu_n) - \phi(s) \right)  + f(s,\phi(s+\mu_n)) \\ && -  f(s,\phi(s)) + g(s+\mu_n) - g(s) \big\| ds  \\
& \geq & H_0 - 2 r R_0 \geq \frac{H_0}{2}.
\end{eqnarray*}
Additionally, the inequality
\begin{eqnarray*} 
\left\|\phi(\tau_n+\mu_n+)-\phi(\tau_n+)\right\| \leq \left(\left\|I+B\right\|+L_h\right) \left\|\phi(\tau_n+\mu_n)-\phi(\tau_n)\right\| 
\end{eqnarray*} 
implies for $t \in [\tau_n-r, \tau_n]$ that
\begin{eqnarray*}
\left\|\phi(t+\mu_n)-\phi(t)\right\| & \geq & \left\|\phi(\tau_n+\mu_n)-\phi(\tau_n)\right\|  \\ 
& - & \displaystyle \Big|\int_{\tau_{n}}^t  \big\| A \left(\phi(s+\mu_n) - \phi(s) \right)  + f(s,\phi(s+\mu_n)) \\ && -  f(s,\phi(s)) + g(s+\mu_n) - g(s) \big\| ds  \Big| \\
& \geq & \frac{H_0}{\left\|I+B\right\|+L_h} -2r R_0 \\
& \geq & \frac{3 H_0}{4\left(\left\|I+B\right\|+L_h\right)}.
\end{eqnarray*} 
For that reason, one can confirm that
\begin{eqnarray*} 
\left\|\phi(t+\mu_n)-\phi(t)\right\| \geq \frac{H_0}{2} \min\left\{1, \frac{3}{2\left(\left\|I+B\right\|+L_h\right)}\right\} \geq \epsilon_0
\end{eqnarray*}
for $t \in [\tau_n-r, \tau_n+r]$.

In the cases I, II, and III, for each natural number $n$, we proved the existence of a real number $\tau_n \geq \theta_{\eta_n p}$ such that (\ref{thmproof2}) holds for $t \in [\tau_n-r, \tau_n+r]$. The sequence $\left\{\tau_n\right\}_{n \in \mathbb N}$ diverges to infinity since the same is true for the sequence $\left\{\eta_n\right\}_{n \in \mathbb N}.$ Consequently, the bounded solution $\phi(t)$ of the impulsive system (\ref{mainimpulsive}) is unpredictable. $\square$
 
\begin{remark}
According to equation (\ref{unpconst}) the unpredictability constant $\epsilon_0$ of $\phi(t)$ is proportional to the unpredictability constant $\delta_0$ of the sequence $\left\{\sigma_k\right\}_{k \in \mathbb N}.$
\end{remark}

\section{An Example} \label{examplessec}

It was proved in paper \cite{Fen17} that the logistic map
\begin{eqnarray} \label{logisticmap1}
z_{k+1}= \Gamma z_k (1-z_k),
\end{eqnarray}
where $k \in \mathbb Z$, possesses an unpredictable orbit for the values of the parameter $\Gamma$ between $3+(2/3)^{1/2}$ and $4$. The unit interval $[0,1]$ is invariant under the iterations of (\ref{logisticmap1}) for these values of the parameter \cite{Hale91}. 

Let us fix an unpredictable orbit $\left\{\sigma^*_k\right\}_{k \in \mathbb Z}$ of the map (\ref{logisticmap1}) with $\Gamma=3.95$, whose terms belong to the interval $[0,1]$. The bounded sequence $\left\{\sigma_k\right\}_{k \in \mathbb Z}$, where $$\sigma_k=(4.5\sigma^*_k, (\sigma^*_k+1)^3), \ k \in \mathbb Z,$$
is unpredictable according to Theorem 3.2 \cite{Fen18}. 

Now, we consider the impulsive system
\begin{eqnarray} \label{examplemain1}
\begin{array}{l}
\displaystyle \frac{dx_1}{dt}= - 6x_1 + 2x_2 + 0.1 \cos(x_2) + 0.3 \sin (2t) + g_1(t), \\
\displaystyle \frac{dx_2}{dt}= - 8x_1 + x_2 + 0.2 \tanh (x_1) + g_2(t), ~t \neq \theta_k, \\
\Delta x_1 |_{t = \theta_k} = - \displaystyle \frac{2}{3} x_1(\theta_k) +0.05 \arctan \left(x_1(\theta_k)\right)+0.4, \\
\Delta x_2 |_{t = \theta_k} = - \displaystyle \frac{2}{3} x_2(\theta_k) +0.04 \sin \left(x_2(\theta_k)\right),
\end{array}
\end{eqnarray}
where $\theta_k=\displaystyle \frac{1}{2} \left((-1)^k+ \pi k\right)$, $k \in \mathbb Z$. In (\ref{examplemain1}), the functions $g_1(t)$ and $g_2(t)$ are defined by the equations $g_1(t)=4.5\sigma^*_k$ and $g_2(t)=(\sigma^*_k+1)^3$ for $t \in (\theta_{2k}, \theta_{2(k+1)}]$, $k \in \mathbb Z$. For each integer $k$, the equation (\ref{regularitycond}) is satisfied for $p=2$ and $\omega=\pi$.

The impulsive system (\ref{examplemain1}) is in the form of (\ref{mainimpulsive}) with
$$A=\begin{pmatrix}  -6 & 2 \\ -8 & 1 \end{pmatrix}, \ \ B=\begin{pmatrix}  -\displaystyle\frac{2}{3} & 0 \\ 0 & -\displaystyle\frac{2}{3}  \end{pmatrix},$$
$$f(t,x_1,x_2)=(0.1 \cos(x_2) + 0.3 \sin (2t), 0.2 \tanh(x_1)), \ \ g(t)=(g_1(t),g_2(t)),$$
and
$$h(x_1,x_2)=(0.05 \arctan \left(x_1 \right)+0.4, 0.04 \sin \left(x_2 \right)).$$
The eigenvalues of the matrix
$$A+\displaystyle \frac{p}{\omega} \log(I+B) = \begin{pmatrix}  -6 - \displaystyle\frac{2}{\pi} \ln 3 & 2 \\ -8 & 1- \displaystyle\frac{2}{\pi} \ln 3 \end{pmatrix}$$
are $\displaystyle - \frac{5}{2}- \frac{2}{\pi} \ln 3 \pm i \frac{\sqrt{15}}{2} $ such that their real parts are negative. Moreover, the equation
$$U(t,s)=e^{-\frac{5}{2} (t-s)} \left(\displaystyle \frac{1}{3}\right)^{i([s,t))} P 
\begin{pmatrix} \displaystyle \cos(\frac{\sqrt{15}}{2}(t-s)) & \displaystyle -\sin( \frac{\sqrt{15}}{2} (t-s)) \\ \displaystyle \sin(\frac{\sqrt{15}}{2}(t-s)) & \displaystyle \cos(\frac{\sqrt{15}}{2}(t-s)) \end{pmatrix} P^{-1}, \ \ t>s,$$
holds, where $U(t,s)$ is the matriciant of the linear homogeneous impulsive system
\begin{eqnarray*} 
&& \displaystyle \frac{dx_1}{dt}= - 6 x_1 + 2 x_2, \\
&& \displaystyle \frac{dx_2}{dt}= - 8 x_1 + x_2, ~t \neq \theta_k, \\
&& \Delta x_1 |_{t = \theta_k} = - \displaystyle \frac{2}{3} x_1(\theta_k), \\
&& \Delta x_2 |_{t = \theta_k} = - \displaystyle \frac{2}{3} x_2(\theta_k) 
\end{eqnarray*}
and $P=\begin{pmatrix}  0 & 1 \\ \displaystyle \frac{\sqrt{15}}{4} & \displaystyle \frac{7}{4} \end{pmatrix}$.

It can be deduced that the assumptions $(A1)-(A7)$ are satisfied for (\ref{examplemain1}) with $N=4.9625$, $\lambda=2.5$, $M_f=0.4473$, $M_h=0.4803$, $L_f=0.2$, and $L_h=0.05$. Therefore, the impulsive system (\ref{examplemain1}) possesses a unique asymptotically stable unpredictable solution according to Theorem \ref{impulsethm}.

\begin{figure}[ht!] 
\centering
\includegraphics[width=13.45cm]{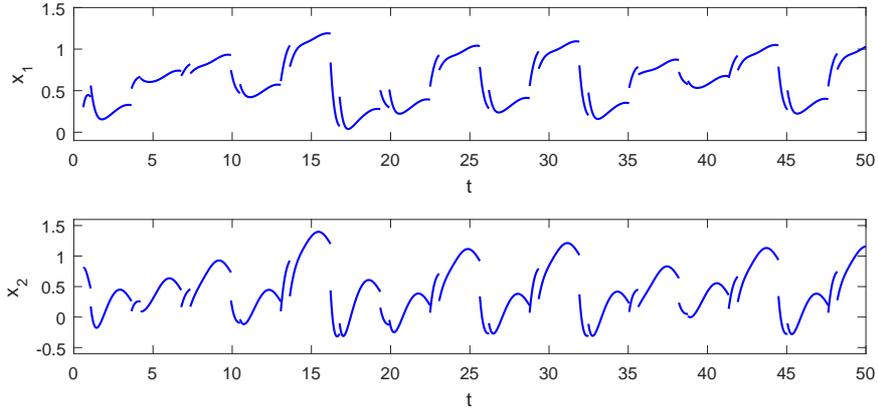}
\caption{Time series of the $x_1$ and $x_2$ coordinates of the solution of (\ref{examplemain1}) in which $g_1(t)=4.5\widetilde{\sigma}_k$ and $g_2(t)=(\widetilde{\sigma}_k+1)^3$ for $t \in (\theta_{2k}, \theta_{2(k+1)}]$, $k \in \mathbb Z$. The figure reveals the existence of a discontinuous unpredictable solution of the impulsive system (\ref{examplemain1}). The initial data $x_1(0.6)=0.3$, $x_2(0.6)=0.8$ are utilized in the simulation.}
\label{impfig}
\end{figure}

In order to demonstrate the unpredictable behavior numerically, we fix an orbit $\left\{\widetilde{\sigma}_k\right\}_{k \in \mathbb Z}$  of the logistic map (\ref{logisticmap1}) with $\Gamma=3.95$ satisfying $\widetilde{\sigma}_0=0.23$, and set $g_1(t)=4.5\widetilde{\sigma}_k$, $g_2(t)=(\widetilde{\sigma}_k+1)^3$ for $t \in (\theta_{2k}, \theta_{2(k+1)}]$, $k \in \mathbb Z$, in  system (\ref{examplemain1}). Figure \ref{impfig} depicts the time series for both coordinates of the solution of (\ref{examplemain1}) corresponding to the initial data $x_1(0.6)=0.3$ and $x_2(0.6)=0.8$. The irregular behavior of the solution confirms the existence of a discontinuous unpredictable solution of (\ref{examplemain1}).

\section{Concluding Remarks} \label{secconc}

In this study, asymptotically stable unpredictable solutions of quasilinear systems with regular moments of impulses is investigated for the first time in the literature. It is worth noting that the presence of nonlinear terms in both the differential equation as well as the impulse equation in system (\ref{mainimpulsive}) is one of the novelties of this study. Another novelty is the usage of regular moments of impulses in the system under consideration. 

The impulsive system (\ref{mainimpulsive2}) admits a unique periodic solution, which is asymptotically stable, provided that the assumptions $(A1)-(A6)$ are fulfilled \cite{Samolienko95}. For that reason, the main source of the unpredictable behavior of system (\ref{mainimpulsive}) is the unpredictable sequence $\left\{\sigma_k\right\}_{k \in \mathbb Z}$, which is used to generate the perturbation $g(t)$. Taking into account the problem under consideration from the input-output mechanism point of view, one can conclude that if an input function of the form $g(t)$ is applied to (\ref{mainimpulsive2}), then an unpredictable output is obtained. In the future, the obtained results can be developed for impulsive systems with variable moments of impulses \cite{Akh1}.


\section*{Conflict of Interest}

This work does not have any conflicts of interest.

\section*{Orcid}
\noindent Mehmet Onur Fen https://orcid.org/0000-0002-7787-7236 \\
\noindent Fatma Tokmak Fen https://orcid.org/0000-0002-4051-7798

\end{document}